\documentclass{amsart}
\usepackage[utf8]{inputenc}
\usepackage[T1]{fontenc}
\usepackage{amsmath}
\usepackage{amssymb}
\usepackage{etaremune}
\usepackage{hyperref}
\usepackage{xcolor}
\usepackage{amsthm}
\usepackage{comment}

\newtheorem{theorem}{Theorem}
\newtheorem{remark}{Remark}

\newtheorem{claim}{Claim}

\newtheorem{lemma}{Lemma}
\newtheorem{definition}{Definition}

\newtheorem{corollary}{Corollary}
\newtheorem{problem}{Problem}

\usepackage{enumitem}

\usepackage{tgbonum}

\usepackage{natbib}
\usepackage{graphicx}

\newcommand{\Cay}{\operatorname{Cay}}
\newcommand{\Sch}{\operatorname{Sch}}
\newcommand{\cost}{\operatorname{Cost}}

\newcommand{\dist}{\operatorname{dist}}
\newcommand{\spn}{\operatorname{Span}}
\newcommand{\subcost}{\operatorname{Sub-Cost}}

\begin{document}

\title{A counterexample to the subgraphings approximate the cost problem}

\author{Bal\'azs Bursics}
	
	\address{E\"otv\"os Lor\'and University, Institute of Mathematics, P\'azm\'any P\'eter stny. 1/C, 1117 Budapest, Hungary}
    \email{bursicsb@gmail.com}

    \author{Jan Greb\'ik}
	
	\address{Charles University, Computer Science Institute, Malostranské nám. 25
118 00 Praha 1, Czechia}
    \email{grebikj@gmail.com }

	\author{Zolt\'an Vidny\'anszky}
    
	\address{E\"otv\"os Lor\'and University, Institute of Mathematics, P\'azm\'any P\'eter stny. 1/C, 1117 Budapest, Hungary}
    \email{zoltan.vidnyanszky@ttk.elte.hu}

\begin{abstract}
    We prove that for every $k\ge 2$, the free group $\mathbb{F}_k$ on $k$ generators admits a finite generating set $S$ such that the Schreier graphing of the Bernoulli action of $\mathbb{F}_k$ with respect to $S$ does not admit spanning subgraphings that would approximate the cost of the action.
    In particular, this gives an example of a graphing that spans a treeable countable pmp equivalence relation, but does not admit a spanning subtreeing.
\end{abstract}

\maketitle

\section{Introduction}

An important numerical invariant of {\bf \emph{probability measure preserving countable Borel equivalence relations (pmp CBERs)}}, introduced by Levitt \cite{Levitt} and subsequently developed into a far-reaching theory by Gaboriau \cite{GaboriouMeasured,GaboriauCost} (see, e.g., \cite{abert2012rank,VoronoiTesselation,HutchcroftPete,le2014number,khezeli2025products,slutsky2026product}), is the notion of {\bf \emph{cost}}.
On a high-level, cost measures how difficult it is to realize a given pmp CBER as connectivity components of a measurable graph relation (also known as graphing) in terms of its average degree.

The main motivation for our work stems from problems about global behavior of cost.
For instance, \emph{how difficult it is to compute a cost of a given pmp action}, or \emph{is it possible to approximate the cost of a countable pmp equivalence relation induced by a graphing by erasing edges}.
In this paper, we describe an example that gives a negative answer to the latter problem, that we will refer to as {\bf \emph{the subgraphings approximate the cost problem}}, and that has been formally stated in \cite[Subgraphings vs Cost Problem, Question~2.68]{GaboriauNotes2025} and \cite[Question~46]{Abert2010}.
Before stating our main result, we formally introduce the main concepts.

\medskip

A {\bf\emph{graphing}} on a standard probability space $(X,\mu)$ is a triplet $\mathcal{G}=(X,E,\mu)$, where $(X,E)$ is a locally countable Borel graph on $X$ that satisfies the mass transport principle with respect to $\mu$, see \cite[Chapter~18]{L12}.
We write $\mathcal{E}(\mathcal{G})$ for the connectedness equivalence relation on $(X,\mu)$ induced by $\mathcal{G}$, it follows that $\mathcal{E}(\mathcal{G})$ is a pmp CBER.
Given a pmp CBER $\mathcal{E}$ on a standard probability space $(X,\mu)$, we define $\spn(\mathcal{E})$ to be the set of all graphings $\mathcal{G}$ on $(X,\mu)$ such that $\mathcal{E}=\mathcal{E}(\mathcal{G})$ up to $\mu$-null sets.
In the special case, when $\mathcal{E}=\mathcal{E}(\mathcal{G})$, we write simply $\spn(\mathcal{G})$ instead of $\spn(\mathcal{E}(\mathcal{G}))$.

A prominent family of graphings are the so-called Schreier graphings of pmp actions: we say that an action $\Gamma\curvearrowright^\sigma (X,\mu)$ of a countable group $\Gamma$ on a standard probability space $(X,\mu)$ is a {\bf \emph{probability measure preserving (pmp) action}} if $\mu(\gamma\cdot A)=\mu(A)$ for every measurable set $A\subseteq X$ and $\gamma\in \Gamma$, where $\gamma\cdot A=\{\gamma\cdot x:x\in A\}$.
Additionally, we say that the action $\sigma$ is (almost everywhere) {\bf \emph{free}} if the stabilizer of $\mu$-almost every point $x\in X$ is trivial.
We denote by $\mathcal{E}(\sigma)$ the orbit equivalence relation induced by $\sigma$, it follows that $\mathcal{E}(\sigma)$ is a pmp CBER.
We set $\spn(\sigma)=\spn(\mathcal{E}(\sigma))$.

For a free pmp action $\Gamma\curvearrowright^\sigma (X,\mu)$ and a finite set $S\subseteq \Gamma\setminus \{1_{\Gamma}\}$, we define
$$\Sch(\sigma,S)=\left(X,\{\{x,\gamma\cdot x\}:x\in X, \ \gamma\in S\},\mu\right)$$
to be the {\bf \emph{Schreier graphing of $\sigma$ with respect to $S$}}.
For this definition, we do not require $S$ to be a generating set.
Observe that $S$ generates $\Gamma$ if and only if $\Sch(\sigma,S)\in \spn(\sigma)$.

\medskip

For a graphing $\mathcal{G}$, we define its {\bf \emph{cost}} as
$$\cost(\mathcal{G})=\frac{1}{2}\int_X \deg_\mathcal{G}(x) \ d\mu.$$
A {\bf \emph{cost of a pmp CBER $\mathcal{E}$}} is then defined as
$$\cost(\mathcal{E})=\inf_{\mathcal{G}\in \spn(\mathcal{E})}\cost(\mathcal{G}),$$
and a {\bf \emph{cost of an action $\sigma$}} is defined as $\cost(\sigma)=\cost(\mathcal{E}(\sigma))$.
Given a graphing $\mathcal{G}$, we define the {\bf \emph{subgraphing cost}}
$$\subcost(\mathcal{G})=\inf\left\{\cost(\mathcal{H}):\mathcal{H}\in \spn(\mathcal{G}) \ \& \ \mathcal{H}\subseteq \mathcal{G}\right\}.$$
Note that we always have $\cost(\mathcal{E}(\mathcal{G}))\le \subcost(\mathcal{G})$, and the subgraphings approximate the cost problem asks whether we have  equality for all graphings $\mathcal{G}$.
If $S$ generate $\Gamma$, then we write simply $\subcost(\sigma,S)$ instead of $\subcost(\Sch(\sigma,S))$.

\medskip

The main focus of this paper is on the free group
$$\mathbb{F}_k=\langle\alpha_1,\dots,\alpha_k\rangle$$
that is freely generated by $\alpha_1,\dots,\alpha_k$ for $k\ge 2$.
It follows from Gaboriau's theorem about treeings \cite[Théorème IV.1]{GaboriauCost} that $\cost(\sigma)=k$ for any free pmp action $\mathbb{F}_k\curvearrowright^\sigma (X,\mu)$.
Our main result is the following.

\begin{theorem}\label{thm:Main}
    Let $k\ge 2$, $\sigma$ be the Bernoulli action $\mathbb{F}_k\curvearrowright \left([0,1]^{\mathbb{F}_k},\lambda^{\otimes\mathbb{F}_k}\right)$ and $\epsilon>0$.
    Then there is a finite set of generators $S_\epsilon$ of size $k+1$ such that
    $$\subcost(\sigma,S_\epsilon)\ge k+1-\epsilon.$$
    In particular, if $\epsilon<1$, then $\subcost(\sigma,S_\epsilon)>\cost(\sigma)=k$.
\end{theorem}

Let $\mathcal{E}$ be a pmp CBER on a standard probability space $(X,\mu)$.
We say that $\mathcal{E}$ is \emph{\bf treeable}, if there is an acyclic graphing $\mathcal{G}$ on $(X,\mu)$ such that $\mathcal{E}=\mathcal{E}(\mathcal{G})$ holds $\mu$ almost everywhere.
The following is a direct consequence of Theorem~\ref{thm:Main} and \cite[Théorème IV.1]{GaboriauCost}.

\begin{corollary}\label{cor:Treeable}
	There exists a graphing $\mathcal{G}$ such that $\mathcal{E}(\mathcal{G})$ is treeable and there is no $\mathcal{H}$ that is an almost everywhere acyclic subgraphing of $\mathcal{G}$ such that $\mathcal{E}(\mathcal{G})=\mathcal{E}(\mathcal{H})$.
\end{corollary}

The following is a direct consequence of the case $k=2$ of Theorem~\ref{thm:Main} (see Section \ref{sec:ProofMain}).

\begin{corollary}\label{cor:MST}
    Assume that $\Lambda$ is a finitely generated group. Let $\Gamma=\mathbb{F}_2*\Lambda$, $\sigma$ be the Bernoulli action $\Gamma\curvearrowright \left([0,1]^\Gamma,\lambda^{\otimes\Gamma}\right)$ and $\epsilon>0$.
    Then there is a graphing $\mathcal{G}_\epsilon\in \spn(\sigma)$ such that
    $$\subcost(\mathcal{G}_\epsilon)\geq\cost(\sigma)+1-\epsilon.$$
\end{corollary}

\begin{remark}
    \label{r:contrast}
    In contrast, recently Recke and the second author have shown \cite[Theorem~1.5]{Poisson} that Schreier graphs of the Bernoulli actions of co-compact lattices in connected higher rank semisimple real Lie groups with property (T) have the \emph{sparse unique infinite cluster property} for any finite set of generators.
    The same conclusion can be derived for Schreier graphs of the Bernoulli actions of products of free groups with at least two factors on the same number of generators $\mathbb{F}_d\times\dots\times\mathbb{F}_d$, for $d\ge2$, by a combination of \cite[Theorem~1]{PoissonProduct} and \cite[Corollary~9.4]{Poisson}.
    It follows that the subgraphings approximate the cost problem has a positive answer for these Schreier graphings. 
\end{remark}

Another application concerns the so-called {\bf \emph{free uniform spanning forest (FUSF)}}.
Let $G=(V,E)$ be an infinite connected locally finite graph, and let $(V_n)_{n\in\mathbb N}$ be an increasing exhaustion of $V$ by finite sets such that each induced subgraph $G_n=G\upharpoonright V_n$ is connected.
Let $\mu_n$ be the law of a uniformly chosen spanning tree of $G_n$, viewed as a probability measure on $\{0,1\}^E$ by declaring all edges outside $G_n$ absent.
The measures $\mu_n$ converge weakly to a probability measure $\mathfrak{F}_G$ on $\{0,1\}^E$.
This limit is independent of the exhaustion and is invariant under every automorphism of $G$; see \cite[Section~5]{BLPS2001}, and it is called the \emph{FUSF of $G$}.

Pete and Tim\'ar \cite[Theorem~1.2]{PeteTimar2022} showed that connectivity of the FUSF can depend on the generating set.
They constructed a virtually free group with two finite generating sets such that the FUSF of one Cayley graph is connected almost surely, whereas that of the other has infinitely many components almost surely.
This answered negatively Question~10.50 of \cite{LyonsPeres2016}.
Combining our main result with the recent theorem that the FUSF is a factor of iid \cite[Theorem~1.1]{openAIFUSF}, we obtain the analogous result for free groups.

\begin{theorem}\label{thm:FUSF}
    Let $k\ge2$.
    Then there are two generating sets $S_0$ and $S_1$ of $\mathbb{F}_k$ such that if we write $G_i=\Cay(\mathbb{F}_k,S_i)$ for $i\in \{0,1\}$, then $\mathfrak{F}_{G_0}$ is connected almost surely and $\mathfrak{F}_{G_1}$ has infinitely many connectivity components almost surely.
\end{theorem}
\begin{proof}
    Let $S_0=\{\alpha_1,\dots,\alpha_k\}$ be the standard generating set for $\mathbb{F}_k$.
    It follows that $\mathfrak{F}_{G_0}$ is concentrated on the standard Cayley graph of $\mathbb{F}_k$, hence it is connected almost surely.

    Let $S_1=S_\epsilon$ be the generating set from Theorem~\ref{thm:Main} for $0<\epsilon<1$.
    By \cite[Theorem~1.1]{openAIFUSF}, $\mathfrak{F}_{G_1}$ is a $\mathbb{F}_k$-invariant factor of iid.
    Equivalently, it can be realized as a measurable acyclic subgraphing $\mathcal T\subseteq\Sch(\sigma,S)$ on the
    Bernoulli probability space.
    As $\mathcal{T}\subseteq \Sch(\sigma,S_1)$, we have that $\mathcal{T}\not\in \spn(\sigma)$ as otherwise we have
    $$\cost(\mathcal T)=\cost(\sigma)=k$$
    by Gaboriau's theorem \cite[Théorème IV.1]{GaboriauCost}.
    In particular, $\mathfrak{F}_{G_1}$ is disconnected with positive probability, hence almost surely as it is ergodic.
    Finally, by \cite[Theorem~1.3]{HutchcroftNachmias2017}, the FUSF of a Cayley graph is almost surely either connected
    or has infinitely many connected components.
\end{proof}

\subsection{Outline of the proof of Theorem~\ref{thm:Main}}\label{subsec:HighLevel}
We only focus on the case $k=2$, as the general case follows directly from Corollary~\ref{cor:MST}.
Write $\mathbb{F}_2=\langle a,b\rangle$.
Given $\epsilon>0$, the generating set $S_\epsilon$ is of the form $\{a,w_\epsilon,bw_\epsilon\}$ for some reduced word $w=w_\epsilon$ with the following combinatorial properties:
\begin{enumerate}
    \item [(A)] large collections of prefixes of $w$ freely generate a free subgroup of $\mathbb{F}_2$ (see, Lemma~\ref{l:freely}),
    \item [(B)] for any fixed $B\subseteq \mathbb{F}_2\setminus \{1_{\mathbb{F}_2}\}$, if we consider the graph $\mathcal{H}_B$ on $\mathbb{F}_2$ with edge set
    $$\{\{\gamma,a\gamma\}:\gamma\in \mathbb{F}_2\}\cup\{\{\gamma,w\gamma\}:\gamma\in \mathbb{F}_2\}\cup\{\{\gamma,b\gamma\}:\gamma\in B\},$$
    then any path from $1_{\mathbb{F}_2}$ to $b$ in $\mathcal{H}_B$ must pass through a $w$-edge $\{\gamma,w\gamma\}$ with the property that almost all (depending on $\varepsilon$) of its prefixes that are followed by letter $b$ belong to $B$ (see, Lemma~\ref{lem:Bypass}).
\end{enumerate}
The word $w$ is constructed by a simple randomized process, namely, it is of the form
$$w=d_{\ell}d_{\ell-1}\dots d_1,$$
where the parameter $\ell\to \infty$ as $\epsilon\to 0$ and $d_i\in \{a,b\}$ are picked uniformly at random.

The above properties are used to derive a result that goes seemingly in the opposite direction than Theorem~\ref{thm:Main}.
To be more precise, we show that if $Z_b\subseteq [0,1]^{\mathbb{F}_2}$ is a measurable set such that
$$\mathcal{H}_{Z_b}=\Sch(\sigma,\{a,w\})\cup\{\{x,b\cdot x\}:x\in Z_b\}\in \spn(\sigma),$$
then $\lambda^{\otimes \mathbb{F}_2}(Z_b)\to 1$ as the parameter $\ell$ defining the length of $w$ goes to infinity (see, Theorem~\ref{thm:Combo}).

On a high-level, the argument uses (B) to gradually enlarge $Z_b$ by ``filling in'' the small amount of missing $b$-edges in selected $w$-edges that eventually produces a positive fraction of $w$-edges whose prefixes are completely contained in $Z_b$ (see, Lemma~\ref{lm:filling}).
If $Z_b$ were small, then this concentration of $Z_b$ on a positive fraction of $w$-edges would contradict an entropy inequality for a restriction of the Bernoulli action to a free subgroup of $\mathbb{F}_2$ (see, Section~\ref{sec:EntropyIneq}).
Here, it is important that the rank of this free group goes to $\infty$ as $\ell\to \infty$, which is where (A) is used.

Finally, Theorem~\ref{thm:Main} is a direct corollary of Theorem~\ref{thm:Combo}.
For example, a lower bound on the measure of $bw$-edges that are used in a spanning subgraphing $\mathcal{G}\subseteq \Sch(\sigma,\{a,w,bw\})$ can be drawn by taking a union with $\Sch(\sigma,\{a,w\})$ and replacing every $bw$-edge $\{x,bw\cdot x\}$ in $\mathcal{G}$ by a $b$-edge $\{w\cdot x,b\cdot (w\cdot x)\}$.
This preserves connectedness, and the measure of $b$-edges remains the same as the measure of $bw$-edges in $\mathcal{G}$.
The lower bound is then derived from Theorem~\ref{thm:Combo}, see Section~\ref{sec:ProofMain}.
A similar argument can be applied to $a$ and $w$.

\subsection{AI Disclosure}
An interesting candidate for a counterexample to the subgraphings approximate the cost problem was suggested by GPT-5.6 in August 2026, upon the authors' inquiry about a construction related to the question of the complexity of finite cost (see Section \ref{sec:prob}).
While neither the authors nor GPT 5.6 were able to prove anything about the group, the newly released GPT 6 showed, interestingly, that the candidate either fails the fixed price conjecture or is a counterexample to the subgraphings approximate the cost problem.\footnote{We still don't know which of the two options holds for that group.}
Then, upon seeing the new OpenAI counterexample to the fixed price conjecture we used extensively the LLM to find a group (it was $\mathbb{Z}*(\mathbb{Z}/3\mathbb{Z})$) that satisfies the above theorem.
After several rounds of drastic simplifications of the arguments, it turned out that the same result is true for the free group, which is what we decided to include.
LLM was not used in the actual writing of this text, only in corrections.

\section*{Acknowledgments:}
We are very grateful to M\'arton Borb\'enyi, Edward Hou, G\'abor Pete for the conversations and suggestions on the topic of the paper.

JG was supported by an ERC CZ grant LL2608 funded by the Czech Ministry of Education, Youth and Sports from the ERC CZ Program. BB and VZ were supported by the ERC COG grant Define (No. 101231345). Views and opinions expressed are however those of the author(s) only and do not necessarily reflect those of the European Union or the European Research Council Executive Agency.

\section{Positive words in $\mathbb{F}_2$}
As discussed in Section~\ref{subsec:HighLevel}, the generating set in Theorem~\ref{thm:Main} is always of the form $\{a,w,bw\}$ for some reduced word $w$; now we describe how to construct them. We use a simple random construction ensuring that overlaps between the word and its translated copies are small.
This implies, among other things, that there is a large collection of prefixes of $w$ that freely generate a free subgroup of $\mathbb{F}_2$.

\begin{definition}[Positive words in $\mathbb{F}_2$]
    Let $\ell$ and $m$ be positive integers.
    An {\bf \emph{$(\ell,m)$-positive word}} is a word $w$ of the form
    \begin{equation}\label{eq:WordForm}
    d_{\ell}d_{\ell-1} \dots d_1,
    \end{equation}
    where $d_i\in \{a,b\}$ for every $1\le i\le \ell$, such that
    \begin{enumerate}
        \item [(a)] the collection
    \begin{equation}\label{eq:DifferentCycles}
        \left\{(d_i,\dots,d_{i+m-1}):1\le i\le \ell\right\}
    \end{equation}
    consists of pairwise distinct sequences, where the indices are always understood $\mod \ell$,
    \item [(b)] we have
    \begin{equation}\label{eq:LLN}
        \ell/3\le |\{i\in \{1,\dots,\ell\}:d_i=a\}|\le 2\ell/3.
    \end{equation}
    \end{enumerate}
    
\end{definition}

Given $d=(d_1,\dots,d_\ell)\in \{a,b\}^\ell$, we write $w_d$ for the reduced word given by \eqref{eq:WordForm}.
We start by showing that a uniform sequence in $\{a,b\}^\ell$ defines a positive word with high probability.
Throughout the paper we use the natural logarithm.

\begin{lemma}
	\label{l:positive}
    Let $\ell$ and $m$ be positive integers such that $5\log(\ell)\le m \le \ell/2$, and endow $\{a,b\}^\ell$ with the uniform measure.
    Then
    $$\mathbb{P}\left(w_d \text{ is an } (\ell,m)\text{-positive word}\right)\ge 1-10/\ell.$$
\end{lemma}
\begin{proof}
    Let $d$ be sampled uniformly from $\{a,b\}^\ell$ and fix $1\le i\not=j\le \ell$.
    The probability that $(d_i,\dots,d_{i+m-1})=(d_j,\dots,d_{j+m-1})$ is equal to $2^{-m}$.
    Indeed, we may without loss of generality assume that $i=1$ and $2\le j\le\lfloor\ell/2\rfloor+1$.
    Then it follows that the choices of $d_1,d_2,\dots,d_{k}\in \{a,b\}$ fully determine $(d_j,\dots,d_{j+m-1})$, where $k=\min\{m,j-1\}$, as $m\le \ell/2$.
    The union bound gives
    $$\mathbb{P}\left(w_d \text{ fails condition (a)}\right)\le 2^{-m}\binom{\ell}{2}\le \ell^{-5\log(2)}\frac{\ell^2}{2}\le \ell^{-3.4}\frac{\ell^2}{2}\le  \frac{1}{\ell}$$
    which gives an upper bound on the probability that (a) fails.

    To bound the probability of failure in (b), we apply Chebyshev's inequality.
    Define the random variable $N_a=|\{i\in \{1,\dots,\ell\}: d_i=a\}|$ and note that $\mathbb{E}[N_a]=\ell/2$ and $\operatorname{Var}(N_a)=\ell/4$.
    Then we have
    $$\mathbb{P}\left(\left|N_a-\frac{\ell}{2}\right|\ge\frac{\ell}{6}\right)\le \frac{9}{\ell},$$
    and the claim follows.   
\end{proof}

Let $w_d$ be an $(\ell,m)$-positive word for some $d\in \{a,b\}^\ell$.
Set $c_1=1_{\mathbb{F}_2}$ and define
$$c_i=d_{i-1}\dots d_1$$
for every $2\le i\le \ell$.
For $u\in \{a,b\}$, let
$$I_u=\left\{i\in \{1,\dots,\ell\}: d_i=u\right\}.$$
The following is the main combinatorial observation.
We state it for the left-shift action of ${\mathbb{F}_2}$ on itself, but later we will apply it to the (orbits of the) Bernoulli shift of $\mathbb{F}_2$.

\medskip

We fix some terminology that is used throughout.
Recall that a \emph{\bf path} in a graph is an injective sequence of vertices, where consecutive vertices form an edge.

\begin{lemma}\label{lem:Bypass}
    Let $\ell$ and $m$ be positive integers, $w_d$ be an $(\ell,m)$-positive word given by some $d\in \{a,b\}^\ell$, and $B\subseteq \mathbb{F}_2\setminus \{1_{\mathbb{F}_2}\}$.
    Define a graph $\mathcal{H}_B$ on $\mathbb{F}_2$ with an edge set
    $$E(\mathcal{H}_B)=\{\{\gamma,a\gamma\}:\gamma\in \mathbb{F}_2\}\cup \{\{\gamma,w_d\gamma\}:\gamma\in \mathbb{F}_2\}\cup \{\{\gamma,b\gamma\}:\gamma\in B\}$$
    and suppose that $1_{\mathbb{F}_2}=\gamma_1,\dots,\gamma_k=b$ is a path in $\mathcal{H}_B$.
    Then there are $1\le i\le k-1$ and $\eta\in \mathbb{F}_2$ such that
    $$(\gamma_i,\gamma_{i+1})=(\eta, w_d\eta) \text{ or } (\gamma_i,\gamma_{i+1})=(w_d\eta,\eta)$$ and 
    \begin{equation}\label{eq:SmallSet}
        1\le \left|\left\{j\in I_b:c_j\eta\not\in B\right\}\right|\le 2m.
    \end{equation}
\end{lemma}
\begin{proof}
    Let $\Cay(\mathbb{F}_2)$ be the left Cayley graph of $\mathbb{F}_2$ with respect to the generating set $\{a,b\}$.
    We start with the following observation.

    \begin{claim}\label{cl:Intersections}
        Let $\eta,\eta'\in \mathbb{F}_2$ and suppose that
        \begin{equation}\label{eq:Intersections}
        \left|\left\{j\in I_b:\exists j'\in I_b \ (c_j\eta,bc_j\eta)=(c_{j'}\eta',bc_{j'}\eta')\right\}\right|>m.
        \end{equation}
        Then $\eta=\eta'$.
    \end{claim}
    \begin{proof}
        Let $j_0\in I_b$ be minimal such that $(c_{j_0}\eta,bc_{j_0}\eta)=(c_{j'_0}\eta',bc_{j'_0}\eta')$ for some $j_0'\in I_b$.
        By \eqref{eq:WordForm}, we have that if $j>j_0$ is such that $j\in I_b$ and $(c_j\eta,bc_j\eta)=(c_{j'}\eta',bc_{j'}\eta')$ for some $j'\in I_b$, then $j'-j'_0=j-j_0$ and
        $$(d_{j_0+1},\dots,d_j)=(d_{j_0'+1},\dots,d_{j'_0+j-j_0}).$$
        By \eqref{eq:Intersections}, there is $j_1\in I_b$ such that $|\{j\in I_b:j_0\le j\le j_1\}|> m$.
        This can only happen, by \eqref{eq:DifferentCycles}, if $j_0=j_0'$.
        Consequently $c_{j_0}\eta=c_{j_0}\eta'$, which implies that $\eta=\eta'$ as desired.    
    \end{proof}
    
    We consider a quotient graph of $\Cay(\mathbb{F}_2)$ that is obtained by gluing together vertices of the form $\{\gamma,a\gamma\}$ for every $\gamma\in \mathbb{F}_2$ and $\{\gamma,b\gamma\}$ for $\gamma\in B$.
    Formally, let $V$ be the set of equivalence classes of the equivalence relation on $\mathbb{F}_2$ that is generated by the relation
    $$\left\{\{\gamma,a\gamma\}:\gamma\in \mathbb{F}_2\right\}\cup\left\{\{\gamma,b\gamma\}:\gamma\in B\right\}.$$
    For $\gamma\in \mathbb{F}_2$, we denote by $[\gamma]$ the corresponding element in $V$ that contains $\gamma$.
    We define the edge relation as
    $$R=\left\{\{[\gamma],[b\gamma]\}:\gamma\not\in B\right\},$$
    and think of each edge $\{v,w\}\in R$ to have a label of the form $\{\gamma,b\gamma\}$, where $\{[\gamma],[b\gamma]\}=\{v,w\}$.

    For $\eta\in \mathbb{F}_2$, define $T_\eta$ to be an increasing enumeration of the set
    $$\{j\in I_b: c_j\eta\not\in B\}.$$
    We denote by $\dist_\mathcal{Q}({-},{-})$ the graph distance in $\mathcal{Q}=(V,R)$.

    \begin{claim}\label{cl:Tree}
        The graph $\mathcal{Q}=(V,R)$ is acyclic, does not contain loops and the label of each edge is uniquely determined.
        Additionally, $\left([\eta],[bc_{j_1}\eta],\dots,[bc_{j_{|T_\eta|}}\eta]\right)$, where $T_\eta=(j_1,\dots,j_{|T_\eta|})$ is a path in $\mathcal{Q}$ that connects $[\eta]$ and $[w_d\eta]$ for every $\eta\in \mathbb{F}_2$.
        In particular, $\dist_\mathcal{Q}([\eta],[w_d\eta])=|T_\eta|$ for every $\eta\in \mathbb{F}_2$.
    \end{claim}
    \begin{proof}
        The first assertion follows from the fact that $\Cay(\mathbb{F}_2)$ is acyclic and $\mathcal{Q}$ is a quotient of $\Cay(\mathbb{F}_2)$ with the property that the restriction of $\Cay(\mathbb{F}_2)$ to each $[\gamma]$, for $\gamma\in \mathbb{F}_2$, is connected.
         The additional parts follow directly from these properties of $\mathcal{Q}$ and the definition.
    \end{proof}

    Given a path $P=(a_1,\dots,a_s)$ in $\mathcal{Q}$, we define $P^r=(a_s,\dots,a_1)$ to be the reverse of $P$, and if $Q=(b_1,\dots,b_t)$ is another path in $\mathcal{Q}$ such that $b_1=a_s$, then we define $Q\circ P$ to be the shortest path from $a_1$ to $b_t$.
    Note that $Q\circ P$ is well defined by Claim~\ref{cl:Tree}, and that it can be equivalently defined as
    \begin{equation}\label{eq:Concat}
        Q\circ P=(a_1,\dots,a_{s-s'(Q, P)},b_{s'(Q, P)+2},\dots,b_t),
    \end{equation}
    where $0\le s'(Q, P)\le s-1$ is maximal such that $b_{j+1}=a_{s-j}$ for every $0\le j\le s'(Q, P)$.
    In particular, $|Q\circ P|=|Q|+|P|-2(s'(Q, P)+1)+1$.

    Suppose for a contradiction that \eqref{eq:SmallSet} does not hold for any $1\le i\le k-1$.
    This is equivalent to $\dist_\mathcal{Q}([\gamma_i],[\gamma_{i+1}])$ either being $0$ or at least $2m+1$ for every $1\le i\le k-1$.
    We prove that there is an $1\le i_0\le k$ such that 
    \begin{equation}    \label{e:obs}
            \dist_{\mathcal{Q}}([1_{\mathbb{F}_2}],[\gamma_i])\ge 2m+1, \text{ for every } i_0\le i\le k.
    \end{equation}
    Once we have that, the proof is finished as $2m+1\le \dist_\mathcal{Q}([1_{\mathbb{F}_2}],[\gamma_k])=\dist_\mathcal{Q}([1_{\mathbb{F}_2}],[b])=1$ contradicts the assumption that $m$ is a positive integer.

    \medskip

    Define
    $$J=\left\{i\in \{1,\dots,k-1\}:[\gamma_i]\not=[\gamma_{i+1}]\right\},$$
    and write $J^+$ for those indices $i\in J$ such that $\gamma_{i+1}=w_d\gamma_i$.
    Set $J^-=J\setminus J^+$ and note that $[\gamma_i]=[w_d\gamma_{i+1}]$ for every $i\in J^-$.
    As $[\gamma_k]=[b]\not=[1_{\mathbb{F}_2}]$, we have that $J\not=\emptyset$.
    We show that setting $i_0=1+\min J$ works as required.
    
    \medskip

    For $i\in J$, define $P_i=(v^i_1,\dots,v^i_{n(i)})$ to be the path in $\mathcal{Q}$ that connects $[\gamma_i]$ and $[\gamma_{i+1}]$ from Claim~\ref{cl:Tree}.
    Note that by our assumption, we have that $|P_i|=n(i)\ge 2m+2$.
    
    For $i_0\le i\le k$, define $J(i)=\max\{j\in J:j<i\}$ and set $Q_i=(w^i_1,\dots,w^i_{m(i)})$ to be the shortest path in $\mathcal{Q}$ from $[\gamma_i]$ to $[1_{\mathbb{F}_2}]$.
    Note that $Q_i$ is uniquely determined by Claim~\ref{cl:Tree}.

    \medskip
    
    We show by induction on $i_0\le i\le k$, that
    \begin{enumerate}
        \item [(a)] $|Q_i|\ge 2m+2$,
        \item [(b)] $w^i_j=v^{J(i)}_{n(J(i))-(j-1)}$ for $1\le j\le m+2$.
    \end{enumerate}
    Note that this will yield \eqref{e:obs}, and hence prove the lemma.
    If $i=i_0$, then (a) and (b) follow as $Q^r_{i_0}=P_{i_0-1}=P_{J(i_0)}$ and $|P_i|\ge 2m+2$ holds for every $i\in J$.
    The inductive step $i\mapsto i+1$ is trivial if $i\not\in J$ as $J(i)=J(i+1)$ and $Q_i=Q_{i+1}$.
    Hence we are left with the case $i\mapsto i+1$, where $i\in J$.

    \begin{claim}\label{cl:Cancel}
        The largest initial segment on which $P_{i}$ and $P^r_{J(i)}$ coincide has size at most $m+1$.
    \end{claim}
    \begin{proof}
        Let $v_0$ and $v_1$ be consecutive vertices in both $P_i$ and $P^r_{J(i)}$.
        Set $\eta=\gamma_{J(i)}$ if $J(i)\in J^+$ and $\eta=\gamma_{J(i)+1}$ otherwise.
        Similarly, set $\eta'=\gamma_i$ if $i\in J^+$ and $\eta'=\gamma_{i+1}$ otherwise.
        
        By the uniqueness of the labels in $\mathcal{Q}$, Claim~\ref{cl:Tree}, and the definition of $P_i$ and $P_{J(i)}$, there are indices $s\in T_{\eta}$ and $t\in T_{\eta'}$ such that
        $$\{c_s\eta,bc_s\eta\}=\{c_t\eta',bc_t\eta'\}.$$
        Note that this implies that
        $$v_r=[c_s\eta]=[c_t\eta'] \text{ and } v_{1-r}=[bc_s\eta]=[bc_t\eta']$$
        for some $r\in \{0,1\}$.
        Indeed, otherwise we get the relation $c^{-1}_tbc_s=(bc_t)^{-1}c_s$ which is not satisfied in $\mathbb{F}_2$.

        Suppose for a contradiction that the initial segment on which $P_{i}$ and $P^r_{J(i)}$ coincide has size strictly more than $m+1$.
        By the previous paragraph and the uniqueness of labels in $\mathcal{Q}$, we conclude that
        $$\left|\left\{j\in I_b:\exists j'\in I_b \ (c_j\eta,bc_j\eta)=(c_{j'}\eta',bc_{j'}\eta')\right\}\right|>m.$$
        Hence, by Claim~\ref{cl:Intersections}, we have $\eta=\eta'$.
        As $\gamma_1,\dots,\gamma_k$ is a path and $J(i)<i$, we must have $J(i)=i-1$, $J(i)\in J^-$, $i\in J^+$ and $\gamma_i=\eta$.
        This implies that $\gamma_{i-1}=\gamma_{i+1}$ and that is a contradiction.
    \end{proof}

    By the definition, we have $Q_{i+1}=Q_i\circ P_{i}^r$, and by the inductive assumption (b), we have that the largest initial segment where $Q_i$ and $P^r_{J(i)}$ coincide has size at least $m+2$.
    By Claim~\ref{cl:Cancel}, we conclude that $s'(Q_i,P_{i}^r)\le m$.
    The description of $Q_i\circ P_{i}^r$ in \eqref{eq:Concat} implies directly that the inductive assumptions (a) and (b) hold for $i+1$, hence finishing the proof of Lemma \ref{lem:Bypass}.
\end{proof}

    Next, we show that a sufficiently spaced collection of initial segments of an $(\ell,m)$-positive word $w_d$ from $I_b$ freely generates a free group.
    Write $I_b=\{t_1<\dots<t_r\}$, where $r=|I_b|$.
    Set $p=\left\lfloor\frac{r}{2(m+1)}\right\rfloor$ and $j_i=t_{2(m+1)i}$ for every $1\le i\le p$.
    Define
	\begin{equation}
	\label{eq:Gdef}
	{\bf F}=\{c_{j_1},\dots,c_{j_p}\}.
	\end{equation}
    Note that by (b), we have that $p\ge \left\lfloor\frac{\ell}{6(m+1)}\right\rfloor$.
	
	\begin{lemma}
		\label{l:freely}
		The elements of the set ${\bf F}$ freely generate a free subgroup of $\mathbb{F}_2$.
	\end{lemma}
	\begin{proof}
        Wlog, we may assume $p \geq 2$. Set $h_1=c_{j_1}$, and for $2\le s\le p$ define
		\begin{equation*}
		\begin{split}
		h_s= & \ c_{j_{s}}(c_{j_{s-1}})^{-1}=d_{j_{s}-1}\dots d_{j_{s-1}}.
		\end{split}
		\end{equation*}
		Observe that each $h_s$ contains at least $2m+1$ occurrences of $b$.
        We say that a $b$ letter is in the middle of the word $h_s$ if there are at least $m$ other $b$ letters to the left and to the right of it.
		We claim that $\{h_1,\dots,h_p\}$ freely generate a free group.
		Note that once we have that, the proof is done as they generate the same subgroup as the elements of ${\bf F}$ that has the same size.
		
		Let $W=h^{\epsilon_1}_{i_1}\cdot \dots h^{\epsilon_n}_{i_n}$ be a reduced word in the alphabet $\{h_1,\dots,h_{p}\}$, where $\epsilon_k\in \{1,-1\}$ for every $1\le k\le n$.
		Suppose for a contradiction that $W=1_{\mathbb{F}_2}$.
		Consider the cancellation process that turns $W$ into $1_{\mathbb{F}_2}$.
		There must be a first time in this process, when one of the $b$ letters in the middle of some of the words $\{h^{\epsilon_k}_{i_k}\}_{k=1}^n$ is canceled.
		Fix such an index $1\le k\le n$.
		As it is the first time in the process, we may assume without loss of generality that the cancellation came from $h^{\epsilon_{k+1}}_{i_{k+1}}$.
        
		Note that we must have $\epsilon_{k+1}\not=\epsilon_k$ as otherwise no cancellation in $h^{\epsilon_{k}}_{i_{k}}h^{\epsilon_{k+1}}_{i_{k+1}}$ occurs at all.
		This implies that $h_{i_k}\not=h_{i_{k+1}}$, and condition (a) in the definition of an $(\ell,m)$-positive word guarantees that at most $m$ letters in $h^{\pm 1}_{i_{k}}$ are canceled in $h^{\pm1}_{i_{k}}h^{\mp1}_{i_{k+1}}$. 
    \end{proof}

\section{Entropy inequality in free groups}\label{sec:EntropyIneq}

Let $k\ge 2$ and $\mathbb{F}_k=\langle\alpha_1,\dots,\alpha_k\rangle$ be a free group on $k$ generators.
We state the entropy inequality (or equivalently the non-negativity of the $f$-invariant) of Bowen \cite{Bowen2010}, see also \cite{BackhauszSzegedy2018}, for the Bernoulli shift $\mathbb{F}_k\curvearrowright \left([0,1]^{\mathbb{F}_k},\lambda^{\otimes \mathbb{F}_k}\right)$.
To that end, recall that the {\bf \emph{entropy of a random variable $R$}} with values in a finite set $A$ is defined as
$$H(R)=-\sum_{a\in A} \mathbb{P}(R=a)\log\left(\mathbb{P}(R=a)\right),$$
where we use the natural logarithm and set $0\log(0)=0$.
Now suppose that $R$ is defined on $\left([0,1]^{\mathbb{F}_k},\lambda^{\otimes \mathbb{F}_k}\right)$, and set $R_i(x)=R(\alpha_i\cdot x)$ for every $x\in [0,1]^{\mathbb{F}_k}$ and $1\le i\le k$.
The {\bf \emph{entropy inequality}} then reads as follows:
\begin{equation}\label{eq:EntropyInequality}
    (2k-1)H(R)\le \sum_{i=1}^k H(R,R_i),
\end{equation}
where $(R,R_i)$ is the joint random variable of $R$ and $R_i$.

\begin{theorem}\label{thm:EntropyIneq}
    Let $\epsilon>0$, $\mathbb{F}_k=\langle\alpha_1,\dots,\alpha_k\rangle$ be a free group on $k\ge 2$ generators and $C\subseteq S\subseteq [0,1]^{\mathbb{F}_k}$ be measurable sets such that $\lambda^{\otimes \mathbb{F}_k}(S)\le 1-\epsilon$ and $\alpha_i\cdot x\in S$ for every $x\in C$ and $1\le i\le k$.
    Then
    $$\lambda^{\otimes \mathbb{F}_k}(C)\le \frac{\log(3)}{k\epsilon}.$$
\end{theorem}
\begin{proof}
    Set $c=\lambda^{\otimes \mathbb{F}_k}(C)$ and $s=\lambda^{\otimes \mathbb{F}_k}(S)$.
    We may assume that $s>0$ as the claim is trivial otherwise.
    Define the random variable $R$ with values in  $\{0,1,2\}$ by
    \[R(x)=\begin{cases}
0& if \ x\in C,\\
1& if \ x\in S \setminus C, \\
2& if \ x \not \in S.
\end{cases}\]
Observe that $H(R)\le \log (3)$.

    We claim that for each fixed $1\le i\le k$, we have that 
    \begin{equation}\label{eq:UpperBoundEntropy}
        H(R,R_i)\le 2H(R)-c\log(1/s).
    \end{equation}
    Indeed, write $p_a=\mathbb{P}(R=a)=\mathbb{P}(R_i=a)$ and $p_{ab}=\mathbb{P}(R=a,R_i=b)$ for $a,b\in \{0,1,2\}$.
    As $p_a=\sum_{b\in \{0,1,2\}}p_{ab}$ for every $a\in \{0,1,2\}$, we have
    \begin{equation*}
        \begin{split}
            2H(R)-H(R,R_i)= & \ -\sum_{a,b\in \{0,1,2\}} p_{ab}\left(\log(p_a)+\log(p_b)\right)+\sum_{a,b\in \{0,1,2\}}p_{ab}\log(p_{ab})\\
            = & \ -2\sum_{a,b\in \{0,1,2\}} p_{ab}\log(p_a)+\sum_{a,b\in \{0,1,2\}}p_{ab}\log(p_{ab})\\
            = & \ \sum_{a,b\in \{0,1,2\}} p_{ab}\log\left(\frac{p_{ab}}{p_ap_b}\right),
        \end{split}
    \end{equation*}
    where in all of the sums, we only consider non-zero probabilities.
    By concavity of $\log$, we have that
    $$\sum_{b\in \{0,1,2\}} p_{ab}\log\left(\frac{p_{ab}}{p_ap_b}\right)\ge \left(\sum_{b\in \{0,1,2\}}p_{ab}\right)\log\left(\frac{\left(\sum_{b\in \{0,1,2\}}p_{ab}\right)}{\left(\sum_{b\in \{0,1,2\}}p_ap_{b}\right)}\right)=p_a\log\left(\frac{p_a}{p_a}\right)=0$$
    for $a\in \{1,2\}$.
    Finally, as $p_{02}=0$, we get from the same inequality that
    $$\sum_{b\in \{0,1\}} p_{0b}\log\left(\frac{p_{0b}}{p_0p_b}\right)\ge p_0\log\left(\frac{p_0}{p_0s}\right)=c\log(1/s)$$
    as $p_0=c$ and $p_0+p_1=s$ by the definition.
    This gives \eqref{eq:UpperBoundEntropy}.

    Summing \eqref{eq:UpperBoundEntropy} for $1\le i\le k$ and using the entropy inequality \eqref{eq:EntropyInequality}, we conclude that 
    $$(2k-1)H(R)\le 2kH(R)-kc\log(1/s),$$
    which is equivalent to
    $$kc\log(1/s)\le H(R)\le \log(3).$$
    The claim follows by noting that $\epsilon\le \log(1/s)$ as $0<s\le 1-\epsilon$.
\end{proof}

\section{Relative cost of $\Sch(\sigma,\{w\})$ with respect to $\{a,b\}$-edges}

We deduce Theorem~\ref{thm:Main} from the following result, which concerns an opposite question.
Namely, suppose that we commit to a graphing $\mathcal{G}\not\in \spn(\sigma)$.
\emph{How many $\{a,b\}$-edges do we need to add to $\mathcal{G}$ to produce an element of $\spn(\sigma)$?}

\medskip

The graphings that we consider are of the form $\Sch(\sigma,\{w\})$, where $w$ is an $(\ell,m)$-positive word.
We show that for certain choices of the parameters $\ell$ and $m$, the measure of the set of $\{a,b\}$-edges that have to be added to turn $\Sch(\sigma,\{w\})$ into an element of $\spn(\sigma)$ goes to $2$ as $\ell\to\infty$.

\begin{theorem}\label{thm:Combo}
    Let $\sigma$ be the Bernoulli action $\mathbb{F}_2\curvearrowright \left([0,1]^{\mathbb{F}_2},\lambda^{\otimes\mathbb{F}_2}\right)$, $\epsilon>0$, and $\ell$ and $m$ be positive integers such that
    \begin{equation}\label{eq:Assumption}
        120(m+1)(2m+1)\log3<2\epsilon^2\ell.
    \end{equation}
    Then for every $(\ell,m)$-positive word $w$ we have that if
    $$\Sch(\sigma,\{w\})\cup \{\{x,a\cdot x\}:x\in Z_a\}\cup \{\{x,b\cdot x\}:x\in Z_b\}\in \spn(\sigma),$$
    for some measurable sets $Z_a,Z_b\subseteq [0,1]^{\mathbb{F}_2}$, then $\lambda^{\otimes\mathbb{F}_2}(Z_r)>1-\epsilon$ for every $r\in \{a,b\}$.
\end{theorem}

We start with an application of Lemma~\ref{lem:Bypass} to the Bernoulli action. Recall the set $\mathbf{F}$ from \eqref{eq:Gdef}.

\begin{lemma}\label{lm:filling}
    Let $\epsilon>0$, $\ell$ and $m$ be positive integers, $w$ be a $(\ell,m)$-positive word and $Z_b$ be a measurable subset of $[0,1]^{{\mathbb{F}_2}}$ such that
	$$\Sch(\sigma,\{a,w\})\cup \{\{x,b\cdot x\}:x\in Z_b\}\in \spn(\sigma)$$
    and $\lambda^{\otimes \mathbb{F}_2}(Z_b)\le 1-\epsilon$.
    Then there is a measurable set $Z_{0,b}\supseteq Z_b$ such that $\lambda^{\otimes \mathbb{F}_2}(Z_{0,b})\le 1-\epsilon/2$ and
    $$\lambda^{\otimes \mathbb{F}_2}\left(\left\{x\in Z_{0,b}:\forall g\in {\bf F} \ g\cdot x\in Z_{0,b}\right\}\right)\ge \frac{2\epsilon}{5(2m+1)}.$$
\end{lemma}
\begin{proof}
    Let $(A_n)_{n\in \mathbb{N}}$ be a sequence of measurable subsets of $[0,1]^{\mathbb{F}_2}$ such that
    \begin{itemize}
        \item [(a)] for every $n\in \mathbb{N}$, we have $\lambda^{\otimes \mathbb{F}_2}(A_n)<\frac{\epsilon}{10(2m+1)}$,
        \item [(b)] for $\lambda^{\otimes {\mathbb{F}_2}}$-almost every $x\in [0,1]^{\mathbb{F}_2}$, we have $|\{n\in \mathbb{N}:x\in A_n\}|=\infty$.
    \end{itemize}

    Consider the following process that starts with $S_0=Z_b$ and inductively creates a sequence of measurable sets $S_n$ for $n\in \mathbb{N}$ such that $S_n\subseteq S_{n+1}$ by the following rule.
    If $x\in A_n$ is such that
    $$1\le \left|\{j\in I_b:c_j\cdot x\not\in S_n\}\right|\le 2m,$$
    then we add $x$ and $c_j\cdot x$ for every $j\in I_b$ to $S_{n+1}$.
    Note that
    \begin{equation}\label{eq:Telescope}
        \lambda^{\otimes \mathbb{F}_2}(S_{n+1})-\lambda^{\otimes \mathbb{F}_2}(S_{n})\le \frac{\epsilon}{10}
    \end{equation}
     for every $n\in \mathbb{N}$ by (a).

    Set $C_n=\left\{x\in S_n:\forall j\in I_b \ c_j\cdot x\in S_n\right\}$.
    Observe that for each $n\in \mathbb{N}$, we have
    $$\lambda^{\otimes \mathbb{F}_2}(S_{n+1})-\lambda^{\otimes \mathbb{F}_2}(S_n)\le (2m+1)\left(\lambda^{\otimes \mathbb{F}_2}(C_{n+1})-\lambda^{\otimes \mathbb{F}_2}(C_n)\right).$$
    Indeed, define a map that sends $y\in S_{n+1}\setminus S_n$ to some $x\in A_n$ such that $$1\le \left|\{j\in I_b:c_j\cdot x\not\in S_n\}\right|\le 2m,$$ and
    $$y\in \{x\}\cup\{c_j\cdot x\}_{j\in I_b}.$$
    This can be done by using the Lusin-Novikov uniformization theorem.
    The claim follows from the fact that such a map is at most $(2m+1)$-to-$1$, and its range falls into $C_{n+1}\setminus C_n$.
    Consequently, we obtain that
    \begin{equation}\label{eq:UB}
        \lambda^{\otimes \mathbb{F}_2}(S_{n})-(1-\epsilon)\le (2m+1)\lambda^{\otimes \mathbb{F}_2}(C_{n})
    \end{equation}
    holds for every $n\in \mathbb{N}$.

    \medskip

    We claim that $\lambda^{\otimes \mathbb{F}_2}(\bigcup_{n\in \mathbb{N}}S_{n})=1$.
    Let $x\in [0,1]^{\mathbb{F}_2}$ be in the free part and such that there is a path $P=(\gamma_1\cdot x,\gamma_2\cdot x,\dots,\gamma_k\cdot x)$, where $\gamma_1=1_{\mathbb{F}_2}$ and $\gamma_k=b$, in
    $$\mathcal{H}_{S_0}=\Sch(\sigma,\{a,w\})\cup\{\{y,b\cdot y\}:y\in S_0\}$$
    from $x$ to $b\cdot x$.
    We claim that $x\in S$, where $S=\bigcup_{n\in \mathbb{N}} S_n$.
    Indeed, suppose for a contradiction that $x\not\in S$.
    As $P$ is also a path in $\mathcal{H}_{S}$ that connects $x$ and $b\cdot x$, we have that there is $y\in P$ such that 
    $$1\le |\{j\in I_b:c_j\cdot y\not\in S\}|\le 2m$$
    by Lemma~\ref{lem:Bypass}.
    Let $p\in \mathbb{N}$ be such that $y\in A_{p}$ and 
    $$\{j\in I_b:c_j\cdot y\not\in S\}=\{j\in I_b:c_j\cdot y\not\in S_{p}\}.$$
    Such an index exists by (b) and the fact that $S=\bigcup_{n\in \mathbb{N}} S_n$.
    This is a contradiction as $\{y\}\cup\{c_j\cdot y\}_{j\in I_b}\subseteq S_{p+1}\subseteq S$ by construction.
    
    \medskip

    By the previous paragraph, there is a maximal index $n_0\in \mathbb{N}$ such that $\lambda^{\otimes \mathbb{F}_2}(S_{n_0})\le 1-\epsilon/2$.
    Set $Z_{0,b}=S_{n_0}$.
    By \eqref{eq:Telescope}, we have
    $$\lambda^{\otimes \mathbb{F}_2}(Z_{0,b})\ge 1-\epsilon/2-\epsilon/10,$$
    which yields, when combined with \eqref{eq:UB},
    $$\frac{2\epsilon}{5(2m+1)}\le \lambda^{\otimes \mathbb{F}_2}(C_{n_0})$$
    as desired.
\end{proof}
    Before turning to the proof of Theorem \ref{thm:Combo} let us make an easy, folklore observation. 

\begin{lemma}
    \label{l:subaction} 
    Let $\Lambda>\Delta$ be countable groups, $\sigma$ be the Bernoulli action $\Lambda\curvearrowright([0,1]^\Lambda,\lambda^{\otimes \Lambda})$, and $\sigma'$ be the $\Delta$ action obtained by restricting $\sigma$ to $\Delta$.
    Then $\sigma'$ is isomorphic to the Bernoulli action of $\Delta$.
\end{lemma}
\begin{proof}
    Choose a collection $R\subseteq\Lambda$ of right coset representatives of $\Delta$.
    The map $\Phi:[0,1]^\Lambda\to ([0,1]^R)^\Delta$ defined by
\[
\Phi(x)(\delta)=(x(\delta r))_{r\in R},
\]
is a $\Delta$-equivariant measure preserving bijection. 
Since $([0,1]^R,\lambda^{\otimes R})$ is a non-atomic standard measure space, it is isomorphic modulo null sets to $([0,1],\lambda)$. Applying this isomorphism in each
$\Delta$-coordinate gives the desired isomorphism.    
\end{proof}

\begin{proof}[Proof of Theorem \ref{thm:Combo}]
    Let $\ell$ and $m$ be positive integers that satisfy
    $$\frac{24(m+1)\log3}{\epsilon\ell}<\frac{2\epsilon}{5(2m+1)},$$
    and $w$ be a $(\ell,m)$-positive word.
    Suppose that there are measurable sets $Z_a,Z_b\subseteq [0,1]^{\mathbb{F}_2}$ such that 
    $$\Sch(\sigma,\{w\})\cup \{\{x,a\cdot x\}:x\in Z_a\}\cup\{\{x,b\cdot x\}:x\in Z_b\}\in \spn(\sigma).$$
    We show that $\lambda^{\otimes \mathbb{F}_2}(Z_r)> 1-\epsilon$ for $r\in \{a,b\}$.
    
    We only handle the case $r=b$, as $r=a$ is completely analogous.
    Suppose for a contradiction that $\lambda^{\otimes \mathbb{F}_2}(Z_r)\le 1-\epsilon$.
    By Lemma~\ref{lm:filling}, there are measurable sets $Z_b\subseteq S$ and $C\subseteq S$ such that $\lambda^{\otimes \mathbb{F}_2}(S)\le 1-\epsilon/2$, $\lambda^{\otimes \mathbb{F}_2}(C)\ge \frac{2\epsilon}{5(2m+1)}$ and $g\cdot x\in S$ for every $g\in {\bf F}$ and $x\in C$.

    \medskip

    Write $\Gamma\simeq\mathbb{F}_k$ for the free group that is freely generated by ${\bf F}$.
    As $\Gamma$ is a subgroup of $\mathbb{F}_2$, we may consider the restriction of $\sigma$ to $\Gamma$, that we denote as $\Gamma\curvearrowright^{\sigma'}[0,1]^{\mathbb{F}_2}$.
    It follows from Lemma \ref{l:subaction} that $\sigma'$ is isomorphic to the Bernoulli shift $\Gamma\curvearrowright [0,1]^\Gamma$ and $C\subseteq S$ satisfy the assumptions of Theorem~\ref{thm:EntropyIneq} with $\epsilon/2$.
    Consequently, $\lambda^{\otimes \mathbb{F}_2}(C)\le \frac{2\log(3)}{k\epsilon}$.
    Finally, note that
    $$k=|{\bf F}|\ge \left\lfloor\frac{\ell}{6(m+1)} \right\rfloor\ge \frac{\ell}{12(m+1)},$$
    where the first inequality is by the definition in \eqref{eq:Gdef} and the second inequality follows from the assumption \eqref{eq:Assumption}.
    This gives the claim.
\end{proof}

\section{Proof of Theorem~\ref{thm:Main}}\label{sec:ProofMain}
In this section we prove our main result.
\begin{proof}[Proof of Theorem \ref{thm:Main}]
The claim is trivial for $\epsilon> 1$ as $\cost(\sigma)=2$, so we may assume that $0<\epsilon\le 1$.
Let $\ell$ and $m$ be positive integers such that
$$\ell\ge 50, \ 5\log(\ell+1)\le m\le \ell/2 \text{ and }    960(m+1)(2m+1)\log3<\epsilon^2\ell.$$
By Lemma~\ref{l:positive}, there is $d\in \{a,b\}^\ell$ such that $w_d$ is an $(\ell,m)$-positive word and $bw_{d}$ is an $(\ell+1,m)$-positive word.
Indeed, the probability that $w_d$ is an $(\ell,m)$-positive word and $bw_d$ is an $(\ell+1,m)$-positive word is at least $1/2-20/\ell>0$ as $\ell\ge 50$. 
Observe that both $w_d$ and $bw_d$ satisfy the assumptions of Theorem~\ref{thm:Combo} for $\epsilon/4>0$.
We show that setting $w=w_d$ and $S_\epsilon=\{a,w,bw\}$ works as required.

The fact that $S_\epsilon$ generates $\mathbb{F}_2$ is clear.
Let $\mathcal{G}\in \spn(\sigma)$ be such that $\mathcal{G}\subseteq \Sch(\sigma,S_\epsilon)$.
Set
$$W_0=\{x\in [0,1]^{\mathbb{F}_2}:\{x,w\cdot x\}\in \mathcal{G}\} \text{ and } W_1=\{x\in [0,1]^{\mathbb{F}_2}:\{x,bw\cdot x\}\in \mathcal{G}\}.$$
Define $Z_1=\{w\cdot x:x\in W_1\}$ and observe that
$$\Sch(\sigma,\{a,w\})\cup \{\{x,b\cdot x\}:x\in Z_1\}\in \spn(\sigma).$$
Indeed, we have that
$$(x,w\cdot x,b\cdot (w\cdot x))$$
is a path in $\Sch(\sigma,\{a,w\})\cup \{\{x,b\cdot x\}:x\in Z_1\}$ for every $x\in W_1$.
By Theorem~\ref{thm:Combo} and the fact that $\lambda^{\otimes\mathbb{F}_2}(Z_1)=\lambda^{\otimes\mathbb{F}_2}(W_1)$, we conclude that $\lambda^{\otimes\mathbb{F}_2}(W_1)>1-\epsilon/4$.
	
Similarly, define $Z_0=\{w\cdot x:x\in W_0\}$.
We claim that
$$\Sch(\sigma,\{a,bw\})\cup \{\{x,b\cdot x\}:x\in Z_0\}\in \spn(\sigma).$$
Indeed, we have that
$$(x,bw\cdot x,w\cdot x)$$
is a path in $\Sch(\sigma,\{a,bw\})\cup \{\{x,b\cdot x\}:x\in Z_0\}$ for every $x\in W_0$.
By Theorem~\ref{thm:Combo} and the fact that $\lambda^{\otimes\mathbb{F}_2}(Z_0)=\lambda^{\otimes\mathbb{F}_2}(W_0)$, we conclude that $\lambda^{\otimes\mathbb{F}_2}(W_0)>1-\epsilon/4$.

Finally, define $A=\{x\in [0,1]^{\mathbb{F}_2}: \{x,a\cdot x\}\in \mathcal{G}\}$, and observe that
$$\Sch(\sigma,\{b,w\})\cup \{\{x,a \cdot x\}:x\in A\}\in \spn(\sigma).$$
By Theorem~\ref{thm:Combo}, we conclude that $\lambda^{\otimes\mathbb{F}_2}(A)>1-\epsilon/2$.

Consequently, we have that
$$\cost(\mathcal{G})\ge \lambda^{\otimes\mathbb{F}_2}(A)+\lambda^{\otimes\mathbb{F}_2}(W_0)+\lambda^{\otimes\mathbb{F}_2}(W_1)>3-\epsilon.$$
As $\mathcal{G}$ was arbitrary, we get that $\subcost(\sigma,S_\epsilon)\ge 3-\epsilon$ as desired.
\end{proof}
Finally, let us deduce Corollary \ref{cor:MST}.

\begin{proof}[Proof of Corollary \ref{cor:MST}]
    Let $c$ be the cost of the equivalence relation given by restricting $\sigma$ to the $\Lambda$ action and let $\mathcal{H}$ be any of its spanning subgraphings. Note that by free-product additivity (see, \cite{GaboriauCost}) $\cost(\sigma)=c+2$. Define $\mathcal{G}_\varepsilon=\mathcal{H} \cup \mathcal{G}'_{\varepsilon}$, where $\mathcal{G}'_{\varepsilon}$ is the graphing provided by Theorem \ref{thm:Main} with $k=2$. Then, since we used a free product, any spanning subgraphing of $\mathcal{G}_\varepsilon$ must span a.e. $\Lambda$ and $\mathbb{F}_2$ component internally. In particular, using Lemma \ref{l:subaction} and our main theorem, its cost is at least $c+2+1-\varepsilon$. 
\end{proof}
\section{Open problems}
\label{sec:prob}
Let us reiterate the question from Kechris-Miller \cite{topics} about the complexity of cost, which was the original motivation behind this paper.

\begin{problem}
Is there a Borel family $(\mathcal{E}_x)_{x \in 2^\mathbb{N}}$ of pmp CBERs such that the set $\{x \in 2^\mathbb{N}: \cost(\mathcal{E}_x)<\infty\}$ is analytic complete?
\end{problem}

The results of this paper and Remark \ref{r:contrast} motivate the following problem.

\begin{problem}
    Characterize the countable groups $\Gamma$, for which the subgraphings approximate cost problem has a positive answer.
    What if we only consider Schreier graphings of the Bernoulli action?
\end{problem}

Note that the answer is positive when $\Gamma$ is amenable: indeed, then its actions are measure hyperfinite. From a hyperfinite decomposition of a graphing $\mathcal{G}$ one can define a spanning subtreeing of $\mathcal{G}$, which in turn must realize the cost of the equivalence relation.

\bibliographystyle{plain}
\bibliography{ref}

\end{document}